\documentclass[11pt]{article}

\usepackage{amsfonts,psfrag}
\usepackage{amsmath, graphicx, epsfig}

\newtheorem{theorem}{Theorem}[section]
\newtheorem{lem}{Lemma}[section]
\newtheorem{cor}{Corollary}[section]

\newcommand{\convas}{\stackrel{a.s.}{\longrightarrow}}

\newcommand{\Nr}{N \rightarrow \infty}

\begin{document}

\title{On the expected time a branching process \\has $K$ individuals alive}
\author{Tom Britton and Peter Neal}
\date{\today}
\maketitle

\begin{abstract} Consider a homogeneous time-continuous
branching process where individuals have constant birth rate
$\delta$, and life length distribution $Q$ having mean $E(Q)=1$. Let
$X(u)$ denote the number of individuals alive at time $u$, and
assume that $X(0)=1$. Let $K$ be a positive integer and define
$A_K:=\int_0^\infty 1_{\{X(u)=K\} }du$, the accumulated time that
the branching process has exactly $K$ individuals alive. In this
paper we prove that $E(A_K)=\delta^{K-1}/\left(
k(1\vee\delta)^K\right)$, irrespective of the life length
distribution $Q$, subject to the normalizing condition $E(Q)=1$.

\end{abstract}


AMS2010 Subject classification: 60J80

\section{Introduction and main result} \label{S:main}

Consider a homogeneous time-continuous branching process
$\mathcal{X}$ having life-length distribution $Q$, where we measure
time in the unit of expected life-lengths, implying that $E(Q)=1$.
During their lives, individuals give birth according to independent
Poisson processes with constant rate $\delta$, where each birth
consists of a single child. Let $X(u)$ denote the number of
individuals alive at time $u$, and assume that the process has one
ancestor, i.e.\ that $X(0)=1$.

Such branching processes have been studied extensively, see, for
example, Jagers \cite{J75} and are sometimes referred to as (binary)
splitting trees, see, for example, Champagnat {\it et al.}
\cite{CLR} . It is well-known that the branching process is
subcritical, critical or supercritical depending on whether $\delta$
is smaller than, equal to, or larger than $1$, and that
$P(X(u)\to\infty)>0$ if and only if $\delta >1$.

Let $N$ denote the number of individuals ever born in the branching
process and let $Q_1, Q_2, \dots $ denote the life-lengths of the
individuals (for example, labelled according to time of birth). It
then holds that
\begin{equation}
\sum_{j=1}^NQ_j=\int_0^\infty X(u)du=\sum_{K=1}^\infty \int_0^\infty
K 1_{\{X(u)=K\} }du ,\label{tot-int}
\end{equation}
where all three expressions equal infinity if the branching process
grows beyond all limits. It is well-known that the expectation of
(\ref{tot-int}) is finite if and only if $\delta<1$, and, using
Wald's lemma,  that the mean then equals
$E(N)E(Q)=E(N)=1/(1-\delta)$.

In the current paper we are interested in the expectation of each of
the terms on the right hand side of (\ref{tot-int}). The factor $K$
is just a constant, so we instead choose to state our result for
\begin{align}
A_K &:= \int_0^\infty 1_{\{X(u)=K\} }du,\quad\text{and its mean} \\
E(A_K)&= E\left( \int_0^\infty 1_{\{X(u)=K\} }du\right) =
\int_0^\infty P(X(u)=K)du,
\end{align}
the expected time the branching process has $K$ living individuals,
or the expected ''$K$-occupation'' time. We have the following
result:

\begin{theorem}\label{th1}
Consider a branching process with birth rate $\delta$ and
life-length distribution $Q$, with $E(Q)=1$ but otherwise arbitrary.
Let $X(u)$ denote the number of individuals alive at time $u$, and
assume $X(0)=1$. Then, for $K=1, 2,\dots$ we have
\begin{equation}
E(A_K) =\int_0^\infty P(X(u)=K)du =\frac{\delta^ {K-1}}{K(1\vee
\delta)^K}.\label{main-result}
\end{equation}
\end{theorem}

The proof of Theorem \ref{th1} is given in Section \ref{S:proof}.

\emph{Remark 1.} The statement of the theorem was conjectured by
Neal \cite{N13} who also proved it for the case that $Q\sim Exp(1)$
and $\delta <1$. In Neal \cite{N13}, it was shown that $E(A_K)$ is
the key quantity for computing the asymptotic endemic equilibrium
distribution for a subcritical branching process with immigration
and a supercritical, homogeneously mixing $SIS$ epidemic model.

\emph{Remark 2.} Theorem \ref{th1} is an example of an {\it
insensitivity} result, in that, $E (A_K)$ only depends upon $Q$
through its mean. Similar results have been observed for queueing
networks, see Zachary \cite{Z07}, for an overview. Furthermore,
Theorem \ref{th1} follows straightforwardly from Zachary \cite{Z07},
Theorem 1 in the subcritical case $\delta <1$. However, Zachary
\cite{Z07}, Theorem 1 does not appear to be easy to adapt to
critical and supercritical branching processes.

The total time until extinction is given by $T=\sum_{k=1}^\infty
A_k =\int_0^\infty 1_{\{X(u)\ge 1\} }du$. It is well-known that this time has infinite mean in the
critical and supercritical case ($\delta \ge 1$). An immediate
consequence of the theorem tells us what the mean equals in the
subcritical case.

\begin{cor} \label{E(T)}
Under the same setting as Theorem \ref{th1} and for the case that
$\delta <1$ we have
\begin{equation}
E(T)=\sum_{K=1}^\infty E(A_K)=\frac{-\log (1-\delta
)}{\delta}.\label{Eq-cor}
\end{equation}
\end{cor}
\emph{Remark}. The case $\delta=0$ can be taken as a limit of
(\ref{Eq-cor}) and the mean is hence $E(T)=1$  as it should.

Theorem \ref{th1} and Corollary \ref{E(T)} have direct implications
for inference on $\delta$ when data consists of observing $A_K$ or
$T=\sum_{K=1}^\infty  A_K$ (cf.\ Farrington and Grant \cite{FG99}).
\begin{cor}\label{est-delta}
Let $K\ge 2$ and suppose that $A_K=t_K$ is observed for the
branching process described above. If the process is known to be
subcritical ($\delta <1$), then the moment estimator of $\delta$ is
given by $\hat\delta=(Kt_K)^{(K-1)^{-1}}$. If the process is super
critical the moment estimator is given by $\hat\delta =1/(Kt_K)$. If
instead $T=t$ is observed and the process is known to be
subcritical, then the moment estimator $\hat\delta$ is the largest
solution to $1-\delta =e^{-\delta t}$.
\end{cor}
\emph{Remark 1.} It might seem unnatural to assume sub- or
super-criticality to be known when making inference. It is however
well-known that branching processes possess many non-standard
inference features (cf.\ Guttorp \cite{G91}).

\emph{Remark 2.} In terms of applications, branching processes are
common  models for populations but also for epidemics, at least when
the outbreaks are small. Corollary \ref{est-delta} treating
estimation problems is hence of interest in these situations, see,
for example, Farrington and Grant \cite{FG99}, treating inference
problems for a related discrete time version of epidemics.

\emph{Remark 3.} The transcendental equation $1-\delta=e^{-\delta
t}$ appears also in mathematical epidemiology, where $t$ is a
measure of infectiousness of the disease and $\delta$ is the
ultimate fraction getting infected in the epidemic. From this theory
(e.g. Diekmann {\it et al.} \cite{DHB12}) it is known that there is
a unique positive solution $\hat\delta $ in (0,1) if and only if
$t>1$.  If $t\le 1$ then the largest solution is $\hat\delta=0$.

Knowing that the expected $K$-occupation time, $E(A_K)$, is
independent of the life-length distribution $Q$, subject to
$E(Q)=1$, raises the question if the result holds even under more
general assumptions. The answer to this question is negative, at
least if we generalize the current model, having constant birth rate
during $Q$, to  a model where the birth rate is inhomogeneous
(time-varying) together with a random duration of the life-length
distribution such that the expected number of births still equals
$\delta$. For this extended model the corresponding $K$-occupation
time is no longer constant (and the same as in Theorem \ref{th1}). A
simple illustration of this is for example $Q\equiv 1$ and $\delta$
large, e.g.\ 10, and $K=1$ and we compare two time-varying birth
rates. The first model is where nearly all of the birth intensity
comes early in life and the second where it comes close to the end
of the life period (i.e. close to 1). Since the process will
probably never return to 1 (being ''very'' super critical), most of
the expected time spent with $K=1$ individual alive comes from
before the first birth, which is clearly longer in the second model.
Also the result no longer holds if we allow the possibility of
multiple births at each point of Poisson point process with rate
$\delta$. It is straightforward to construct a counter example
comparing $Q \equiv 1$ and $Q \sim {\rm Exp} (1)$ since for $Q
\equiv 1$ all individuals born at the same time will die at the same
time.

\section{Proof of Theorem \ref{th1}} \label{S:proof}

\subsection{Introduction}

The approach we take to prove Theorem \ref{th1} is to consider
phase-type distributions (Asmussen {\it et al.}
\cite{Asmussenetal.1996} ) for $Q$. Specifically, we take $Q$ to be
a finite mixture of hypoexponential distributions. That is, we
assume that there exists $m \in \mathbb{N}$ such that $Q$ is a
mixture distribution of $Q_1, Q_2, \ldots, Q_m$ with $P (Q = Q_i)
=p_i$ $(\sum_{i=1}^m p_i =1)$ and, for each $i$, there exists $n_i
\in \mathbb{N}$ such that
\begin{eqnarray} \label{eq:qdist} Q_i \sim {\rm Exp} (\gamma_{i,1}) +  {\rm Exp} (\gamma_{i,2}) + \ldots +  {\rm Exp} (\gamma_{i,n_i}) \end{eqnarray} with
$\sum_{i=1}^m p_i \sum_{j=1}^{n_i} \gamma_{i,j}^{-1}=1$. Therefore
an individual has a lifetime distributed according to $Q_i$ with
probability $p_i$, where the lifetime $Q_i$ consists of $n_i$ stages
each of which lasts an exponential length of time. Thus if we know
the distribution from which each of the individual lifetimes come
and the stage at which each individual is at in their lifetime, the
 branching process is Markovian. Any Coxian distribution
(Asmussen {\it et al.} \cite{Asmussenetal.1996}) can be expressed as
a mixture of hypoexponential distributions, and therefore the above
class of mixtures of hypoexponential distributions is dense,
Asmussen {\it et al.} \cite{Asmussenetal.1996}. Hence Theorem
\ref{th1} follows trivially if we can show that $E (A_K) =
\delta^{K-1}/\{ K (1 \wedge \delta)^K \}$ holds for the mixtures of
hypoexponential distributions.

For $i=1,2,\ldots, m$, let $k_{i,j}$ denote the total number of
individuals with a lifetime distributed according to $Q_i$ who are
currently in stage $j$ of their lifetime with $\mathbf{k}_i =
(k_{i,1}, k_{i,2}, \ldots, k_{i,n_i})$ and $\mathbf{k} =(
\mathbf{k}_1, \mathbf{k}_2, \ldots, \mathbf{k}_m)$. To consider the
transitions to and from state $\mathbf{k}$ it is helpful to define
$\mathbf{e}_{i,j}$ to be a vector of length $\sum_{i=1}^m n_i$ whose
$(i,j)^{th}$ ($\sum_{l=1}^{i-1} n_l +j^{th}$) element is 1 and all
other elements are 0. Transitions from state $\mathbf{k}$ occur as
follows:-
\begin{eqnarray} \label{eq:trans} \mathbf{k} \rightarrow \left\{ \begin{array}{lll} \mathbf{k} - \mathbf{e}_{i,j} +\mathbf{e}_{i,j+1} & \mbox{rate } k_{i,j} \gamma_{i,j} & (i=1,2,\ldots,m; j=1,2, \ldots, n_i -1) \\  \mathbf{k} - \mathbf{e}_{i,n_i} & \mbox{rate } k_{i,n_i} \gamma_{i,n_i} & (i=1,2,\ldots,m) \\
\mathbf{k} + \mathbf{e}_{i,1} & \mbox{rate } K \delta p_i &
(i=1,2,\ldots,m), \end{array} \right. \end{eqnarray} with $K =
\sum_{i=1}^m \sum_{j=1}^{n_i} k_{i,j}$. The Markov branching
(birth-death) process defined by \eqref{eq:trans} will  almost
surely go extinct in the subcritical and critical case. For the
supercritical branching process, the branching process will either
go extinct or grow exponentially large. In all three cases
(subcritical, supercritical and critical branching processes) we
study a modified process which results in an irreducible, aperiodic
and positive recurrent Markov chain. The stationary distribution of
the modified process is the key tool for computing $E (A_K)$. The
details of the appropriate modified process differs between the
three classes of branching process and we therefore consider each in
turn. We use detailed balance to prove the stationary distribution
of the modified process which is that same approach that is often
used in queueing theory to prove insensitivity results, see, for
example, Whittle \cite{Whittle} and Zachary \cite{Z07} and
references therein.

\subsection{Subcritical branching process: $\delta <1$} \label{S:d<1}

Theorem \ref{th1} can be proved in the subcritical case using
Zachary \cite{Z07}, Theorem 1. However, it is instructive for
studying the supercritical case to outline a proof of the result.
The branching process will almost surely go extinct and we modify
\eqref{eq:trans} by regenerating the branching process whenever it
goes extinct by restarting the branching process with a new
individual whose lifetime is distributed according to $Q_l$ with
probability $p_l$. This is the regeneration approach introduced in
Hern\'andez-Su\'arez and Castillo-Chavez \cite{HSCC99} and extended
in Ball and Lyne \cite{BL02}. That is, for $\mathbf{k} =
\mathbf{e}_{i,n_i}$ replace $\mathbf{k} \rightarrow \mathbf{k} -
\mathbf{e}_{i,n_i}$ by $\mathbf{k} \rightarrow \mathbf{e}_{l,1}$
with probability $p_l$ $(l=1,2,\ldots, m)$. The modified process is
an irreducible, aperiodic and positive recurrent Markov chain and
therefore has a unique stationary distribution. Let $\pi_\mathbf{k}$
denote the stationary probability of the population being in state
$\mathbf{k}$.

\begin{lem} \label{lemmab1} For all $\mathbf{k}$,
\begin{eqnarray} \label{eq:stat1}\pi_\mathbf{k} = C (K-1)! \delta^{K-1} \prod_{i=1}^m \prod_{j=1}^{n_i} \frac{q_{i,j}^{k_{i,j}}}{k_{i,j}!},
\end{eqnarray} where  $q_{i,j} = p_i/\gamma_{i,j}$, $C =  -\delta/\log (1-\delta)$ and $K= \sum_{i=1}^m \sum_{j=1}^{n_i} k_{i,j}$.

For $K \geq 1$, let $\mathcal{A}_K = \{ \mathbf{k}: \sum_{i=1}^m
\sum_{j=1}^{n_i} k_{i,j} = K \}$ and
\begin{eqnarray} \label{eq:stat1a} \phi_K = \sum_{\mathbf{k} \in \mathcal{A}_K} \pi_\mathbf{k}  = C \delta^{K-1}/K. \end{eqnarray}

\end{lem}
{\bf Proof.} To prove the Lemma it is sufficient to show that
$\pi_\mathbf{k}$ given by \eqref{eq:stat1} satisfies detailed
balance. That is,
\begin{eqnarray} \label{eq:stat2} \sum_{\mathbf{h} \neq \mathbf{k}} \pi_\mathbf{h} \rho_{\mathbf{h}, \mathbf{k}} = \pi_\mathbf{k} \sum_{\mathbf{l} \neq \mathbf{k}}  \rho_{\mathbf{k}, \mathbf{l}},
\end{eqnarray} where $\rho_{\mathbf{h}, \mathbf{l}}$ is the rate at which transitions from state $\mathbf{h}$ to state $\mathbf{l}$ take place.

We start with the right hand side of \eqref{eq:stat2} in the case,
where $\mathbf{k} \neq \mathbf{e}_{l,1}$ and we do not need to
consider the regeneration modifications.

The possible transitions to state $\mathbf{k}$ are from  $\mathbf{k}
+ \mathbf{e}_{i,j} - \mathbf{e}_{i,j+1}$ $(i=1,2,\ldots, m;
j=1,2,\ldots, n_i-1)$, $\mathbf{k} + \mathbf{e}_{i,n_i}$ and
$\mathbf{k} - \mathbf{e}_{i,1}$. We have the following three
equalities. For $i=1,2,\ldots,m$ and $j=1,2, \ldots,n_i-1$,
\begin{eqnarray} \label{eq:qq1} \pi_{\mathbf{k} + \mathbf{e}_{i,j} - \mathbf{e}_{i,j+1}} (k_{i,j} +1) \gamma_{i,j} &=& \pi_\mathbf{k} \frac{q_{i,j}}{k_{i,j}+1} \frac{k_{i,j+1}}{q_{i,j+1}}
 (k_{i,j} +1) \gamma_{i,j} \nonumber \\
&=&  \pi_\mathbf{k} k_{i,j+1} \gamma_{i,j+1}. \nonumber \\
\end{eqnarray}
 For $i=1,2,\ldots,m$,
\begin{eqnarray} \label{eq:qq2} \pi_{\mathbf{k} + \mathbf{e}_{i,n_i}} (k_{i,n_i} +1) \gamma_{i,n_i} &=& \pi_\mathbf{k} \delta K \frac{q_{i,n_i}}{k_{i,n_i}+1}
 (k_{i,n_i} +1) \gamma_{i,n_i} \nonumber \\
&=&  \pi_\mathbf{k} \delta K p_i, \end{eqnarray} and
\begin{eqnarray} \label{eq:qq3} \pi_{\mathbf{k} - \mathbf{e}_{i,1}}(K-1) p_i \delta &=& \pi_\mathbf{k} \frac{1}{(K-1) \delta} \frac{k_{i,1}}{q_{i,1}} (K-1) p_i \delta \nonumber \\
&=&  \pi_\mathbf{k} k_{i,1} \gamma_{i,1}. \end{eqnarray} It follows
from (\ref{eq:qq1}--\ref{eq:qq3}) that
\begin{eqnarray} \label{eq:stat3} & & \sum_{\mathbf{h} \neq \mathbf{k}} \pi_\mathbf{h} \rho_{\mathbf{h}, \mathbf{k}}  \nonumber \\
&=& \sum_{i=1}^m \sum_{j=1}^{n_i-1} \pi_{\mathbf{k} + \mathbf{e}_{i,j} - \mathbf{e}_{i,j+1}}  (k_{i,j}+1) \gamma_{i,j} +  \sum_{i=1}^m \pi_{\mathbf{k} + \mathbf{e}_{i,n_i}}  (k_{i,n_i}+1) \gamma_{i,n_i} +  \sum_{i=1}^m \pi_{\mathbf{k} - \mathbf{e}_{i,1}} K p_i \delta \nonumber \\
&=& \sum_{i=1}^m \sum_{j=1}^{n_i-1} \pi_{\mathbf{k} }  k_{i,j+1} \gamma_{i,j+1} +  \sum_{i=1}^m \pi_{\mathbf{k}} \delta K p_i +  \sum_{i=1}^m  \pi_\mathbf{k} k_{i,1} \gamma_{i,1} \nonumber \\
&=& \pi_\mathbf{k} \left\{ \sum_{i=1}^m \sum_{j=1}^{n_i} k_{i,j}
\gamma_{i,j} + \sum_{i=1}^m p_i \delta K  \right\} = \pi_\mathbf{k}
\sum_{\mathbf{l} \neq \mathbf{k}}  \rho_{\mathbf{k}, \mathbf{l}} \nonumber \\
\end{eqnarray}
as required.

For $l=1,2,\ldots,m$ and $\mathbf{k} = \mathbf{e}_{l,1}$,  we have
that
\begin{eqnarray} \label{eq:stat4}\sum_{\mathbf{h} \neq \mathbf{e}_{l,1}} \pi_\mathbf{h} \rho_{\mathbf{h}, \mathbf{e}_{l,1}}
&=& \sum_{i=1}^m  \left\{\pi_{\mathbf{e}_{l,1} + \mathbf{e}_{i,n_i}}  \gamma_{i,n_i} + \pi_{\mathbf{e}_{i,n_i}} p_l \gamma_{i,n_i} \right\}    \nonumber \\
&=& \sum_{i=1}^m \left\{ C \delta q_{l,1} q_{i,n_i} \gamma_{i,n_i} + C q_{i,n_i} \gamma_{i,n_i} p_l \right\} \nonumber \\
&=& C (\delta  q_{l,1} + p_l) \sum_{i=1}^m q_{i,n_i} \gamma_{i,n_i} = C (\delta q_{l,1} + 1) \sum_{i=1}^mp_i \nonumber \\
&=& C (\delta q_{l,1} + p_l)  = C q_{l,1} (\delta + \gamma_{l,1} ) \nonumber \\
&=& \pi_{\mathbf{e}_{l,1}}   \sum_{\mathbf{h} \neq \mathbf{e}_{l,1}}
\rho_{\mathbf{e}_{l,1}, \mathbf{h}}
\end{eqnarray} as required.

Since $ \sum_{\mathbf{k} \in \mathcal{A}_K} K! \prod_{i=1}^m
\prod_{j=1}^{n_i} q_{i,j}^{k_{i,j}}/k_{i,j}! = (\sum_{i=1}^m
\sum_{j=1}^{n_i} q_{i,j})^K =1$, \ref{eq:stat1a} follows trivially
from \ref{eq:stat1} and $C =  -\delta/\log (1-\delta)$ follows from
$\sum_{K=1}^\infty \phi_K = C \sum_{K=1}^\infty \delta^{K-1}/K =1$.
 \hfill
$\square$

The above tells us that in equilibrium the modified process spends a
proportion $\phi_K$ of its time with $K$ individuals alive. It does
not tell us directly the mean amount of time, $E (A_K)$, that the
branching process has $K$ individuals. Each cycle, between
regenerations, of the modified process corresponds to a single
realisation of the branching process. As noted in that the mean
total size of the subcritical branching process is $1/(1-\delta)$.
Thus the mean regeneration time, $E(T)$, satisfies
\[ \sum_{K=1}^\infty E(T) \phi_K K = \frac{1}{1 -\delta},\] which gives $E(T) = 1/C= - \log (1- \delta)/\delta$ {\it cf.} Corollary \ref{E(T)}.
Hence $E (A_K) = E(T) \phi_K= \delta^{K-1}/K$ is the mean time the
original subcritical branching process with 1 initial ancestor
spends with $K$ individuals alive.

\subsection{Supercritical case: $\delta > 1$} \label{S:d>1}

We modify the approach taken in Section \ref{S:d<1} to obtain $E
(A_K) = 1/(K \delta)$ when $\delta >1$. In this case the expected
size of the branching process is infinite and there is a  non-zero
probability of never going extinct. This leads to key differences
from the subcritical case which need to be resolved. Therefore we
modify the branching (birth-death) process, defined by
\eqref{eq:trans}, to create a population process, $\mathcal{P}_N$,
which is restricted to $\{1, 2, \ldots, N \}$ individuals for some
$N \in \mathbb{N}$. The population process evolves as the branching
process except for transitions which lead to 0 or $N+1$ individuals
in the branching process. In particular, any transition in the
branching process which leads to $N+1$ individuals (a birth when
there are $N$ individuals) leads to the population process
restarting with 1 individual whose lifetime is distributed according
to $Q$. (A catastrophic event killing all $N$ individuals in the
population combined with a regeneration event, the birth of a new
individual.) Also extinction of the branching process (death of the
only individual) is replaced by the population process moving to a
state with $N$ individuals alive. (This move does not have a natural
interpretation.) The lifetime stages of the $N$ individuals will be
discussed below. The population process is rather different to the
regeneration process used for the subcritical case in Section
\ref{S:d<1} but the resulting Markov chain has a finite state space
and is aperiodic and irreducible  and therefore has a unique
stationary distribution.

\begin{lem} \label{lemmab2} For any $N \in \mathbb{N}$, the population process $\mathcal{P}_N$ with stationary distribution $\{ \pi_\mathbf{k}^N \}$ satisfies the following.

For all $\mathbf{k}$, with $1 \leq K= \sum_{i=1}^m \sum_{j=1}^{n_i}
k_{i,j} \leq N$,
\begin{eqnarray} \label{eq:star1}\pi_\mathbf{k}^N = L_N (K-1)! \prod_{i=1}^m \prod_{j=1}^{n_i} \frac{w_{i,j}^{k_{i,j}}}{k_{i,j}!},
\end{eqnarray} where $L_N =  \left\{ \sum_{j=1}^N (1/j) \right\}^{-1} \approx 1/\log N$ and $\mathbf{w}$ satisfies $\sum_{i=1}^m \sum_{j=1}^{n_i} w_{i,j} =1$,
\begin{eqnarray}
\label{eq:star3a} \gamma_{i,1} + \delta &=& \frac{p_i
\delta}{w_{i,1}} + \sum_{l=1}^m
w_{l,n_l} \gamma_{l,n_l} \\
\label{eq:star3b} \gamma_{i,j} + \delta &=& \frac{\gamma_{i,j-1}
w_{i,j-1}}{w_{i,j}} + \sum_{l=1}^m w_{l,n_l} \gamma_{l,n_l}
\hspace{0.2cm} (j=2,3,\ldots,n_i).
\end{eqnarray}

For $1 \leq K \leq N$, let $\mathcal{A}_K = \{ \mathbf{k}:
\sum_{i=1}^m \sum_{j=1}^{n_i} k_{i,j} = K \}$ then
\begin{eqnarray} \label{eq:star1a} \phi_K^N = \sum_{\mathbf{k} \in \mathcal{A}_K} \pi_\mathbf{k}  = L_N/K. \end{eqnarray}
\end{lem}

We comment briefly on the statement of Lemma \ref{lemmab2} before
embarking on the proof. Firstly, the form of $\pi_\mathbf{k}^N$
given by \eqref{eq:star1} is very similar to that for
$\pi_\mathbf{k}$ given by \eqref{eq:stat1}. The key difference is
that $\delta$ does not feature explicitly in  \eqref{eq:star1} and
$w_{i,j}$ has a more complicated form than $q_{i,j} =
p_i/\gamma_{i,j}$. For the subcritical case $q_{i,j}$ is the
probability that an alive individual is in the $j^{th}$ stage of a
lifetime distributed according to $Q_i$ and $w_{i,j}$ plays the same
role for the supercritical case. In the supercritical case, the
branching process is {\it growing} and we observe that $w_{i,1} >
q_{i,1}$ and $w_{i,n_i} < q_{i,n_i}$. Fortunately it is not
necessary to compute $\mathbf{w}$, although for $Q \sim {\rm Gamma}
(2,2)$ it is straightforward to show that  $w_{1,1} =w$ and $w_{1,2}
= 1-w$, where $w= (- \delta + \sqrt{\delta^2 + 8 \delta})/4$.
Secondly, \eqref{eq:star1a} follows trivially from \eqref{eq:star1}
and gives us that $E(A_K) \propto 1/K$. Therefore after proving
Lemma \ref{lemmab2}, we show that the constant of proportionality is
$1/\delta$.

{\bf Proof of Lemma \ref{lemmab2}.} Fix $N \in \mathbb{N}$. We
follow Lemma \ref{lemmab1} in proving the lemma using detailed
balance. That is, showing that for all $\mathbf{k}$,
$\pi_\mathbf{k}$ given by \eqref{eq:star1} satisfies
\begin{eqnarray} \label{eq:star2} \sum_{\mathbf{h} \neq \mathbf{k}} \pi_\mathbf{h}^N \rho_{\mathbf{h}, \mathbf{k}}
= \pi_\mathbf{k}^N \sum_{\mathbf{l} \neq \mathbf{k}}
\rho_{\mathbf{k}, \mathbf{l}}.
\end{eqnarray}

We start with the right hand side of \eqref{eq:star2} in the case,
where $\mathbf{k} \neq \mathbf{e}_{l,1}$ or $K=N$, that is, we do
not need to worry about modifications.

The possible transitions to state $\mathbf{k}$ are from  $\mathbf{k}
+ \mathbf{e}_{i,j} - \mathbf{e}_{i,j+1}$ $(i=1,2,\ldots, m;
j=1,2,\ldots, n_i-1)$, $\mathbf{k} + \mathbf{e}_{i,n_i}$ and
$\mathbf{k} - \mathbf{e}_{i,1}$. We have the following three
equalities which are similar to  (\ref{eq:qq1}--\ref{eq:qq3}) in
Lemma \ref{lemmab1}. However, there are differences due to the
different form of $\pi_\mathbf{k}$ between the subcritical and
supercritical case. For $i=1,2,\ldots,m$ and $j=1,2, \ldots,n_i-1$,
\begin{eqnarray} \label{eq:qr1} \pi_{\mathbf{k} + \mathbf{e}_{i,j} - \mathbf{e}_{i,j+1}}^N (k_{i,j} +1)
 \gamma_{i,j} &=& \pi_\mathbf{k}^N \frac{w_{i,j}}{k_{i,j}+1} \frac{k_{i,j+1}}{w_{i,j+1}}
 (k_{i,j} +1) \gamma_{i,j} \nonumber \\
&=&  \pi_\mathbf{k}^N k_{i,j+1} \frac{\gamma_{i,j}
w_{i,j}}{w_{i,j+1}}.
\end{eqnarray}
 For $i=1,2,\ldots,m$,
\begin{eqnarray} \label{eq:qr2} \pi_{\mathbf{k} + \mathbf{e}_{i,n_i}}^N (k_{i,n_i} +1) \gamma_{i,n_i} &=&
\pi_\mathbf{k}^N K \frac{w_{i,n_i}}{k_{i,n_i}+1}
 (k_{i,n_i} +1) \gamma_{i,n_i} \nonumber \\
&=&  \pi_\mathbf{k}^N K w_{i,n_i} \gamma_{i,n_i}  \nonumber \\
&=&  \pi_\mathbf{k}^N \sum_{i=1}^m \sum_{j=1}^{n_i} k_{i,j}
w_{i,n_i} \gamma_{i,n_i},
\end{eqnarray} and
\begin{eqnarray} \label{eq:qr3} \pi_{\mathbf{k} - \mathbf{e}_{i,1}}^N(K-1) p_i \delta
&=& \pi_\mathbf{k}^N \frac{1}{(K-1)} \frac{k_{i,1}}{w_{i,1}} (K-1) p_i \delta \nonumber \\
&=&  \pi_\mathbf{k}^N \frac{k_{i,1} p_i \delta}{w_{i,1}}.
\end{eqnarray}

The rate of transition out of state $\mathbf{k}$ is $\sum_{i=1}^m
\sum_{j=1}^{n_i} k_{i,j} \gamma_{i,j} + K \delta =\sum_{i=1}^m
\sum_{j=1}^{n_i} k_{i,j} (\gamma_{i,j} + \delta)$. Hence, for
\eqref{eq:star2} to hold, it follows from
(\ref{eq:qr1}--\ref{eq:qr3}) that we require that
\begin{eqnarray} \label{eq:star3} & & \sum_{i=1}^m
\sum_{j=1}^{n_i-1} k_{i,j+1} \frac{\gamma_{i,j} w_{i,j}}{w_{i,j+1}}
+ \left( \sum_{i=1}^m \sum_{j=1}^{n_i} k_{i,j} \right) \sum_{i=1}^m
w_{i,n_i} \gamma_{i,n_i} + \sum_{i=1}^m \frac{k_{i,1} p_i
\delta}{w_{i,1}} \nonumber \\ &=& \sum_{i=1}^m \sum_{j=1}^{n_i}
k_{i,j} (\gamma_{i,j} + \delta).
\end{eqnarray}
Equating the coefficients of the $k_{i,j}$ terms we get $\mathbf{w}$
solving \eqref{eq:star3a} and \eqref{eq:star3b}. Therefore we need
to check that this choice of $\mathbf{w}$ is also consistent with
the boundary cases.

For $\mathbf{k}$ such that $K=N$, we have that
\begin{eqnarray}
\label{eq:star4} \sum_{l=1}^m \pi_{\mathbf{e}_{l,n_l}}^N
\gamma_{l,n_l} \{ \pi_\mathbf{k}^N/(L_N/N) \} &=& \sum_{l=1}^m L_N
w_{l,n_l} \gamma_{l,n_l} L_N (N-1)! \prod_{i=1}^m \prod_{j=1}^{n_i}
\frac{w_{i,j}}{k_{i,j}!} \times \frac{N}{L_N} \nonumber
\\
&=& \sum_{l=1}^m (k_l+1) \gamma_{l,n_l} \pi_{\mathbf{k} +
\mathbf{e}_{l,n_l}},
\end{eqnarray}
where, with an abuse of notation, we take $\pi_{\mathbf{k} +
\mathbf{e}_{l,n_l}}^N$ to satisfy \eqref{eq:star1} with $K=N+1$.
That is, the transitions from a single individual to $N$ individuals
in the modified process mimic the transitions in the branching
process from $N+1$ individuals to $N$ individuals. All other
transitions in and out of state $\mathbf{k}$ are identical to the
branching process and it is straightforward to verify that
\eqref{eq:star2} holds.

For $\mathbf{k} = \mathbf{e_{l,1}}$, we have that the transitions
into state $\mathbf{e_{l,1}}$ are from $\mathbf{e}_{l,1}+
\mathbf{e}_{i,n_i}$ and from any $\mathbf{k}$ such that $K =N$.
Using $\phi_N^N =\sum_{\mathbf{k} \in \mathcal{A}_N}
\pi_\mathbf{k}^N = L_N /N$, the rate of entry into state
$\mathbf{e}_{l,1}$ is
\begin{eqnarray}
\label{eq:star5} p_l N \delta \phi_N^N + \sum_{i=1}^m \gamma_{i,n_i}
\pi_{\mathbf{e}_{l,1}+ \mathbf{e}_{i,n_i}}^N &=& p_l L_N \delta +
\sum_{i=1}^m \gamma_{i,n_i} \frac{L_N \delta w_{l,1}
w_{i,n_i}}{\delta}. \end{eqnarray} The rate of exit from state
$\mathbf{e}_{l,1}$ is $(\gamma_{l,1} + \delta)$. Hence, from
\eqref{eq:star5} we require that
\begin{eqnarray}
\label{eq:star6} p_l L_N \delta + \sum_{i=1}^m \gamma_{i,n_i} L_N
w_{l,1} w_{i,n_i}= (\gamma_{l,1} + \delta) L_N w_{l,1}.
\end{eqnarray} However, \eqref{eq:star6} is equivalent to
\eqref{eq:star3a}. (Simply multiply both sides of \eqref{eq:star3a}
by $L_N w_{l,1}$.) Therefore completing the proof that
$\pi_\mathbf{k}^N$ is indeed the stationary distribution of the
modified population process.
 \hfill
$\square$

We now show that $E (A_K) = 1/(\delta K)$ by  studying the total
amount of time, $T_N$, that a supercritical branching process spends
with between 1 and $N$ individuals alive. In particular we study the
asymptotic behaviour of $T_N$ as $N \rightarrow \infty$. The
supercritical branching process either goes extinct with
probability, $z$, or explodes ($X (u) \rightarrow \infty$ as $u
\rightarrow \infty$) with probability, $1 -z$. It is helpful to
consider these two cases separately with ${\rm Ext}$ denoting the
event that the branching process goes extinction. Firstly, it
follows from Champagnat {\it et al.} \cite{CLR}, Proposition 2.1,
that $z =1- \eta/\delta$ , (see also Lambert \cite{Lambert}, Section
5), where $\eta$ is the Malthusian parameter of the branching
process. Moreover, Lambert \cite{Lambert}, Proposition 5.7 gives the
law of the supercritical process conditional upon extinction which
is a subcritical branching process with a modified lifetime
distribution and birth rate $\delta - \eta$. Thus it is trivial to
show that $E (T_N | {\rm Ext}) = O(1)$. Conditional upon
non-extinction, Champagnat {\it et al.} \cite{CLR}, Proposition 2.2
states,
 $\exp (- \eta u) X(u) \convas Y$ as $u \rightarrow \infty$, where $Y$ is exponential random variable with mean dependent upon $\eta$ and $Q$.
(This is proved in Lambert \cite{Lambert} and is a special case of
Nerman \cite{Nerman}, Theorem 5.4.)  Let $T_N^F = \min \{u; X(u)
=N+1 \}$ and $T_N^L = \max \{u; X(u) =N \}$ denote the first time
$X(u)$ leaves $\{1,2,\ldots, N \}$ and the last time $X(u)$ belongs
to $\{1,2,\ldots, N\}$ with $T_N^F \leq T_N \leq T_N^L$. Then
conditional upon non-extinction, it is straightforward to show that
$T_N^F \times \eta/\log N, T_N^L \times\eta/\log N | {\rm Ext}^C
\convas 1$ as $\Nr$ and consequently that $T_N \times \eta/\log N (=
\tilde{T}_N) | {\rm Ext}^C \convas 1$ as $\Nr$. Given the structure
of the branching process it is straightforward to show that there
exists $0 < q <1$, such that for all $m \in \mathbb{N}$, $P
(\tilde{T}_N > m) \leq q^m$, and hence, that $E (\tilde{T}_N | {\rm
Ext}^C) \rightarrow 1$ as $\Nr$. Therefore
\begin{eqnarray}
E (T_N) &=& P ({\rm Ext}) E (T_N | {\rm Ext}) +  P ({\rm Ext}^C) E (T_N | {\rm Ext}^C) \nonumber \\ &=& z O(1) + (1-z) \left\{ \frac{\log N}{\eta} + o (\log N) \right\} \nonumber \\
&=&  \frac{\eta}{\delta} \times \frac{\log N}{\eta} + o (\log N) \nonumber \\
&=& \frac{1}{\delta} \log N + o (\log N).
\end{eqnarray}
Since the above holds for all $N$, we have that
\begin{eqnarray}
E(A_K) = \lim_{\Nr} \phi_K^N E (T_N) = \lim_{\Nr}\frac{L_N}{K}
\frac{1}{\delta} \log N = \frac{1}{\delta K}
\end{eqnarray}
as required.

\subsection{Critical case: $\delta = 1$} \label{S:d=1}

For the critical case the expected duration of the branching process
is infinite and the probability of non-extinction is 0, which makes
studying this case particularly difficult. The modified population
process defined in Section \ref{S:d>1} can still be constructed in
this case and Lemma \ref{lemmab2} still holds with $q_{i,j}
(=p_i/\gamma_{i,j}) = w_{i,j}$. Thus it is straightforward to show
that $E (A_K) \propto 1/K$ with  $\lim_{\delta \uparrow 1}
\delta^{K-1}/K = \lim_{\delta \downarrow 1} 1/(\delta K) = 1/K$.
However, it is difficult to adapt the approach taken in Section
\ref{S:d>1} to show that  $E (A_K) =1/K$ but to prove the result it
suffices to fix $Q$ and to show that $E(A_1)=1$.

We study the time that the branching process spends with one
individual alive, which we term the local time process (of the
branching process with one individual). Suppose that $Q$ satisfies
\eqref{eq:qdist}. We say that the local time process is in state
$(i,j)$ if the individual is in the $j^{th}$ stage of lifetime
$Q_i$. The local time process regenerates with a new individual in
state $(i,1)$ with probability $p_i$ $(i=1,2,\ldots,m)$ if the
branching process goes extinct. Also since the branching process is
critical it will almost surely go extinct. Therefore the branching
process started from one individual will eventually return to one
individual, although the mean waiting time is infinite. The
transitions from state $(i,j)$ in the local time process are given
by:-
\begin{eqnarray} \label{eq:cc}
(i,j) \rightarrow \left\{ \begin{array}{ll} (a,b) & \mbox{rate }
r_{(i,j),(a,b)} \\
(i,j+1) & \mbox{rate } \gamma_{i,j} \mbox{ if } (j=1,2,\ldots,n_i-1)
\\
(a,1) & \mbox{rate } \gamma_{i,j} p_a \mbox{ if } j=n_i
\end{array} \right.
\end{eqnarray}
where $\sum_{(a,b)} r_{(i,j),(a,b)} =1$ and $r_{(i,j),(a,b)}$ is the
probability that following a birth to an individual in state $(i,j)$
the next time the branching process returns to one individual, the
sole individual will be in state $(a,b)$. The matrix $R =
(r_{(i,j),(a,b)})$ is difficult to compute but we know the
stationary distribution for the local time process from the
stationary distribution of the branching process. Hence the
probability of being in state $(i,j)$ is $q_{i,j} =
p_i/\gamma_{i,j}$. Therefore following Ball and Lyne \cite{BL02},
Section 2.2.3, we note that the overall intensity of regeneration of
the local time process (extinction of the branching process) is
$\sum_{i=1}^m q_{i,n_i} \gamma_{i,n_i} = \sum_{i=1}^m p_i =1$. Thus
the mean time between regenerations is $E(A_1) =1^{-1} =1$ as
required.

\section*{Acknowledgements}

The paper was written while T.B. spent a sabbatical semester at
University of Florida. T.B. is grateful to Ira Longini and his group
for their hospitality. We would like to thank Frank Ball for
bringing Zachary \cite{Z07} to our attention.


\end{document}